\documentclass[10pt]{amsart}

\usepackage{latexsym,amssymb} 
\usepackage{amsmath,amsfonts,amsthm}
\input xypic

\newtheorem{theorem}{Theorem}[section]
\newtheorem{proposition}[theorem]{Proposition}

\newtheorem{lemma}[theorem]{Lemma}
\newtheorem{example}[theorem]{Example}

\newtheorem{remark}[theorem]{Remark}

\newcommand{\dom}{\mathbf{d}}

\newcommand{\ran}{\mathbf{r}}

\thanks{This research was supported by an EPSRC grant (EP/F004184, EP/F014945, EP/F005881),
the Funda\c{c}\~{a}o para a Ci\^{e}ncia e a Tecnologia, courtesy of Pedro Resende under the grant PPCDT/MAT/55958/2004, 
{\em Groupoids and quantales in geometry and analysis}, at the Instituto Superior T\'ecnico, Lisbon, 
and by Prof Stuart Margolis of Bar-Ilan University, Israel.}

\title{A non-commutative generalization of Stone duality}
\author{M.~V.~Lawson}
\address{Department of Mathematics and the
Maxwell Institute for Mathematical Sciences\\
Heriot-Watt University\\
Riccarton\\
Edinburgh~EH14~4AS\\
Scotland\\
\texttt{M.V.Lawson@ma.hw.ac.uk}}

\begin{document}
\maketitle

\begin{abstract} We prove that the category of boolean inverse monoids is dually equivalent
to the category of boolean groupoids. This generalizes the classical Stone duality between boolean algebras and boolean spaces.
As an instance of this duality, we show that the boolean inverse monoid $C_{n}$ associated with the Cuntz groupoid $G_{n}$ is
the strong orthogonal completion of the polycyclic (or Cuntz) monoid $P_{n}$.
The group of units of $C_{n}$ is the Thompson group $V_{n,1}$.\\

\noindent
2000 {\em Mathematics Subject Classification}: 20M18,18B40,06E15.
\end{abstract}

\section{Statement of the theorem}\setcounter{theorem}{0}

The importance of partial, as opposed to global, symmetries in mathematics is well-established.
The question is how to describe them mathematically.
One approach, advocated in \cite{L}, is to use inverse semigroups; 
these are direct generalizations of groups and are ultimately descended from the pseuodgroups of transformations used in differential geometry.
Another approach is to use topological groupoids such as in the recent work of Hughes \cite{H1,H2}.
Although the surface structure of these two approaches looks very different, they are in fact closely related.
Classically, pseudogroups of tranformations give rise to topological groupoids of germs.
More generally, Paterson \cite{Paterson} used ideas from functional-analysis to construct topological groupoids from inverse semigroups
and Renault \cite{Renault} constructed inverse semigroups from topological groupoids using bisections.
This work has been developed in a number of directions \cite{E1,E2,K1,K2,Lenz,Paterson,R1,R2,S}, to name but a few.
The goal of our paper is to set up an {\em exact} correspondence between a class of inverse monoids, we call boolean monoids, 
and a class of topological groupoids, we call boolean groupoids.
As the terminology suggests, our correspondence can be seen as a natural generalization of the Stone duality between boolean algebras and boolean spaces.
For background on inverse semigroups see \cite{L}, for groupoids \cite{H} and for topological groupoids \cite{E2,R1,Paterson,Renault}.

Although our theorem appears to link semigroups and groupoids in reality it is linking two different kinds of groupoid.
Boolean inverse monoids are semigroups but they are also special kinds of ordered groupoids by virtue of the Ehresmann-Schein-Nambooripad Theorem \cite{L}.
It follows that our theorem could also be viewed as providing a duality between a class of
{\em ordered} groupoids on the one hand and a class of {\em topological} groupoids on the other.
The advantages of being able to re-encode algebraic structures as topological ones cannot be overstated.
In the remainder of the section, we define the two categories we shall work with and put the results of this paper in context.

Let $S$ be an inverse monoid with zero.
If $s \in S$ we write $\dom (s) = s^{-1}s$ and $\ran (s) = ss^{-1}$.
We say that $s,t \in S$ are {\em compatible} if $s^{-1}t$ and $st^{-1}$ are both idempotents and {\em orthogonal} if $s^{-1}t$ and $st^{-1}$ are both zero.
Inverse semigroups come equipped with a partial order, called the {\em natural partial order}, defined by $s \leq t$ iff $s = te$ for some idempotent $e$.
The natural partial order is the only order used in this paper.
If $s$ and $t$ are bounded above they are compatible.
It follows that when discussing the existence of joins in an inverse semigroup we are only interested in elements which are {\em a priori} pairwise compatible.
With respect to the natural partial order, the set of idempotents, $E(S)$, becomes a meet semilattice.
Our perspective is that inverse monoids may therefore be regarded as generalizations of meet semilattices.
The particular inverse monoids considered in this paper have an even stronger order-theoretic character.
We say that an inverse monoid is a {\em boolean inverse monoid} if it satisfies the following three conditions:
\begin{description}

\item[{\rm (BM1)}] $(E(S),\leq)$ is a boolean algebra.

\item[{\rm (BM2)}] $(S,\leq)$ is a meet semilattice.

\item[{\rm (BM3)}] The join of pairs of orthogonal elements always exists.

\end{description}
An inverse semigroup with zero is said to be {\em (finitely) orthogonally complete} if it has 
joins of all finite orthogonal subsets and multiplication distributes over finite orthogonal joins \cite{L2}.
The semilattice of idempotents in a boolean inverse monoid is distributive,
and so using the same argument as Proposition~1.4.20 of \cite{L}
it follows that boolean inverse monoids are orthogonally complete. 
We shall see later, in Lemma~2.3, that in fact such monoids have the joins of all finite non-empty subsets of pairwise compatible elements.

In order to define the morphisms between boolean inverse monoids we need some definitions.
For $A \subseteq S$, define $A^{\uparrow} = \{s \in S \colon \exists a \in A, a \leq s \}$.
If $A = A^{\uparrow}$ we say that $A$ is {\em upwardly closed}.
A {\em filter base} in $S$ is a subset $X \subseteq S$ with the property that $x,y \in X$ implies there exists $z \in X$ such that $z \leq x,y$.
A {\em filter} in $S$ is a subset $F$ which is upwardly closed and a filter base.
In a boolean inverse monoid a filter is an upwardly closed subset closed under finite meets.
A {\em proper filter} is a filter that does not contain zero.
An {\em ultrafilter} in $S$ is a proper filter $F$ which is maximal amongst proper filters. 
If $s \in S, s \neq 0$ then $s^{\uparrow}$ is the {\em principal filter} generated by $s$.

A {\em morphism} $\theta \colon S \rightarrow T$ between two boolean monoids is a semigroup homomorphism such that 
\begin{description}

\item[{\rm (M1)}] $\theta \mid E(S) \colon E(S) \rightarrow E(T)$ is a homomorphism of boolean algebras.
\item[{\rm (M2)}] $\theta \colon (S,\leq) \rightarrow (T,\leq)$ is a homomorphism of semilattices.
\item[{\rm (M3)}]  The inverse image under $\theta$ of every ultrafilter in $T$ is an ultrafilter in $S$.

\end{description}

We now turn to groupoids.
Let $G$ be a groupoid; that is, a small category in which every arrow is invertible.
The set of identities is denoted by $G_{0}$ and the domain and range maps by $\mathbf{d}$ and $\mathbf{r}$, respectively.
We denote by $G \ast G$ the set of composable pairs $(g,h)$ where $\dom (g) = \ran (h)$.
Let $\mathsf{P}(G)$ denote the powerset of $G$.
If $A,B \subseteq G$ define $AB = \{ab \colon a \in A, b \in B, \exists ab \}$ and $A^{-1} = \{a^{-1} \colon a \in A \}$.
With respect to these operations  $\mathsf{P}(G)$ is a semigroup with involution.
An element $A$ of $\mathsf{P}(G)$ is called a {\em bisection} if $a,b \in A$ and $\mathbf{d}(a) = \mathbf{d}(b)$
(respectively $\ran (a) = \ran (b)$) implies that $a = b$.
Equivalently, $A$ is a bisection iff $A^{-1}A,AA^{-1} \subseteq G_{o}$.
We say that a groupoid $G$ is a {\em boolean groupoid} if it satisfies the following conditions:
\begin{description}

\item[{\rm (BG1)}] $G$ is a hausdorff \'etale topological groupoid, where a topological groupoid is {\em \'etale} if its domain map is a local homeomorphism.

\item[{\rm (BG2)}]  $G_{0}$ is compact.

\item[{\rm (BG3)}] $G$ has a basis of compact open bisections.

\end{description}
A {\em morphism} between boolean groupoids is a continuous covering functor.

The main theorem proved in this paper can now be stated.\\

\noindent
{\bf Theorem }{\em The category of boolean inverse monoids is dually equivalent to the category of boolean groupoids}.\\

The motivation for our definition of boolean inverse monoids came from a number of sources.

The pioneering paper on topological groupoids and their connection with inverse semigroups is Renault's \cite{Renault}.
On page~142 \cite{Renault}, he remarks that the inverse monoids constructed from ample topological groupoids
have a boolean algebra of idempotents and finite orthogonal joins.
The importance of the existence of finite orthogonal joins was re-iterated in Paterson's book, see Proposition~4.4.3 \cite{Paterson},
where boolean algebras, in fact generalized boolean algebras, play an important role.
The significance of boolean algebras in the theory of $C^{\ast}$-algebras of topological groupoids has been taken up recently by \cite{E1,E2}.
It was in Proposition~2.9 of \cite{S} that an explicit equivalence was proved between an ample groupoid's being hausdorff
and the existence in the associated inverse semigroup of finite meets; such semigroups were first studied in detail by Leech \cite{Leech}.

The construction of a topological groupoid from an inverse semigroup was first carried out by Renault \cite{Renault}.
It was Paterson \cite{Paterson} who developed Renault's work into a theory of the universal groupoid associated with an inverse semigroup.
Whereas Paterson constructed his groupoid from a functional-analytic perspective,
Lenz \cite{Lenz}, combining ideas from both Paterson and Kellendonk \cite{K1,K2},
showed that one could construct the universal groupoid of an inverse semigroup directly from the inverse semigroup
by using equivalence classes of down-directed subsets of the inverse semigroup.
His motivation seems to have had two sources: first, Kellendonk's technique of building a groupoid using equivalence classes of
descending chains of elements, described in Section~9.2 of \cite{L},
and second, the role played by ultrafilters in the theory of convergence in topological spaces.
Stuart Margolis and the author realized, during a visit of the latter to Bar-Ilan University in January 2009, 
that the equivalence classes Lenz worked with could be replaced by filters.
Thus with each inverse semigroup $S$ one can associate the inverse semigroup of filters $\mathsf{L}(S)$.
By taking the underlying groupoid of this inverse semigroup and introducing a topology derived from the way
$S$ is embedded in $\mathsf{L}(S)$ one gets Paterson's universal groupoid.

Both Paterson and Lenz also constructed a reduction of this groupoid which can, in the language of filters,
be seen as the groupoid of ultrafilters on $S$.
In fact this was what Kellendonk was interested in \cite{K1,K2}.
The importance of the ultrafilters has been taken up in the recent work of Exel \cite{E1,E2}.
Lenz also investigated conditions on the inverse semigroup $S$ which guarantee that this reduced groupoid had pleasant properties.

The idea of trying to prove a duality type theorem linking inverse semigroups and topological groupoids 
arose from conversations with Pedro Resende during the author's visit to Lisbon over Easter 2008.

The catalyst which led to the formulation of the theorem of this paper was \cite{S} which made us realize
that everything in Lenz's paper \cite{Lenz} would work much more easily if only the inverse semigroup had sufficiently rich 
order-theoretic properties.

\section{The proof}\setcounter{theorem}{0}

Our proof is a direct generalization of the familiar proof of the classical Stone duality; see \cite{BS}, for example.
We begin with the algebraic ingredients of our proof.
First, we establish some consequences of the axioms for boolean inverse monoids.
The natural partial order plays a key role and certain properties proved in Section~1.4 of \cite{L} are summarised here.

\begin{lemma} Let $S$ be an inverse semigroup.
\begin{enumerate}

\item $s$ and $t$ are compatible if and only if $s \wedge t$ exists and $\dom (s \wedge t) = \dom (s) \wedge \dom (t)$
and $\ran (s \wedge t) = \ran (s) \wedge \ran (t)$.

\item If $s \vee t$ exists then $\dom (s \vee t) = \dom (s) \vee \dom (t)$ and $\ran (s \vee t) = \ran (s) \vee \ran (t)$.

\item If $s \wedge t$ exists then for any $u \in S$ we have that 
$us \wedge ut$ (respectively $su \wedge tu$) exists and $u(s \wedge t) = us \wedge ut$ (respectively $(s \wedge t)u = su \wedge tu$).

\end{enumerate}
\end{lemma}

Contrast (1) and (2) above: if $s$ and $t$ are not compatible we will {\em not} have both $\dom (s \wedge t) = \dom (s) \wedge \dom (t)$
and $\ran (s \wedge t) = \ran (s) \wedge \ran (t)$.

In an inverse semigroup $S$ we use the notation $s^{\downarrow}$ to mean the set of all elements in $S$ beneath $s$.

\begin{lemma} Let $S$ be a boolean inverse monoid.
\begin{enumerate}

\item For each $s \in S$ the  poset $(s^{\downarrow},\leq)$ is a boolean algebra.

\item Let $s \leq t$.
Then there is a unique element $t \setminus s$ satisfying the following conditions:
$t \setminus s \leq t$, the pair $s$ and $t \setminus s$ are orthogonal and $t = s \vee (t \setminus s)$.

\item  Let $s \neq 0$ and $s \nleq t$. 
Then there exists a non-zero element $s'$ such that $s' \leq s$ and $s' \wedge t = 0$.

\end{enumerate}
\end{lemma}
\proof (1) Define the function $\beta \colon s^{\downarrow} \rightarrow \dom (s)^{\downarrow}$ by $x \mapsto \dom (x)$.
This is an order isomorphism.
Clearly $(\dom (s)^{\downarrow},\leq)$ is a boolean algebra since $(E(S),\leq)$ is.
Thus $(s^{\downarrow},\leq)$ is a boolean algebra.

(2) Put $e = \dom (t) \wedge \dom (s)'$, working in the boolean algebra $E(S)$.
Define $t \setminus s = te$.
By construction $t \setminus s \leq t$ and $\dom (t \setminus s) = e$.
It follows that $s$ and $t \setminus s$ are orthogonal.
Thus their join $s \vee (t \setminus s)$ exists.
Observe that $\dom (s \vee (t \setminus s)) = \dom (t)$ and  clearly $s \vee (t \setminus s) \leq t$.
It follows that $t = s \vee (t \setminus s)$.

Let $x$ be any element such that $x \leq t$, $s$ and $x$ are orthogonal, and $t = s \vee x$.
To show that it is equal to $t \setminus s$ it is enough to show that
$\dom (x) = \dom (t \setminus s)$.
But this follows by the uniqueness of relative complements in boolean algebras.

(3) Let $s,t \in S$ be non-zero elements of a boolean inverse monoid.
Then $s \wedge t \leq s$ and so we may form the element $s' = s \setminus s \wedge t$.
It follows from (2) above that $s' = 0$ iff $s \leq t$.
We deduce that if $s \nleq t$ then there exists a non-zero element $s'$ such that $s' \leq s$ and $s' \wedge t = 0$.\qed \\

As a result of the above lemma, we say that boolean inverse monoids are {\em locally boolean}.
We can now prove that boolean inverse monoids have all finite compatible joins.

\begin{lemma} Let $S$ be a boolean inverse monoid.
If $s$ and $t$ are compatible then $s \vee t$ exists.
\end{lemma}
\proof Both elements $s \setminus s \wedge t$ and $t \setminus s \wedge t$ exist
and so $s = (s \wedge t) \vee (s \setminus s \wedge t)$ and $t = (s \wedge t) \vee (t \setminus s \wedge t)$ by Lemma~2.2.
The elements $s \wedge t$ and $s \setminus s \wedge t$, 
as well as $s \wedge t$ and $t \setminus s \wedge t$ are pairwise orthogonal.
We prove that that $s \setminus s \wedge t$ and $t \setminus s \wedge t$ are orthogonal.
We use the fact that since $s$ and $t$ are compatible, we may apply Lemma~2.1(1).
Thus
$\dom (s \setminus s \wedge t) = \dom (s) \wedge \dom (t)'$ 
and 
$\dom (t \setminus s \wedge t) = \dom (t) \wedge \dom (s)'$;
and
$\ran (s \setminus s \wedge t) = s \dom (t)' s^{-1}$ 
and 
$\ran (t \setminus s \wedge t) = t \dom (s)' t^{-1}$.
It is now clear that the elements are orthogonal.
Put
$$x =  (s \wedge t) \vee (s \setminus s \wedge t) \vee (t \setminus s \wedge t).$$
We prove that $x = s \vee t$.
Clearly $s,t \leq x$ and $\dom (x) = \dom (s) \vee \dom (t)$.
It is easy to check that $x = s \vee t$.\qed \\

We now turn to the properties of filters and ultrafilters on boolean inverse monoids.

If $F$ is a filter in $S$ and $s \in F$ we write $s \wedge F \neq 0$ to mean that $s \wedge a \neq 0$ for all $a \in F$.
The following is a special case of Lemma~12.3 of \cite{E2}.

\begin{lemma} Let $F$ be a filter in the boolean inverse monoid $S$.
Then $F$ is an ultrafilter iff $s \wedge F \neq 0$ implies that $s \in F$.
\end{lemma}

\begin{lemma} Let $S$ be a boolean inverse monoid.
\begin{enumerate}

\item Each non-zero element of $S$ belongs to an ultrafilter.

\item If $a \in S$ is non-zero then the intersection of all ultrafilters containing $a$ is the 
principal filter $a^{\uparrow}$.

\end{enumerate}
\end{lemma}
\proof (1) This is a standard argument using Zorn's Lemma.

(2) Let $A$ be the intersection of all ultrafilters containing $a$.
Clearly $a^{\uparrow} \subseteq A$.
Let $b \in A$.
We prove that $a \leq b$.
Suppose not.
Then by Lemma~2.2(3) there exists an element $a' \neq 0$ such that $a' \leq a$ and $a' \wedge b = 0$.
Let $C$ be any ultrafilter containing $a'$.
Then $a \in C$ and so $A \subseteq C$.
But $a,a',b \in A$ imply that $a',b \in C$ and so $0 = a' \wedge b \in C$, which is a contradiction.
Thus $a \leq b$, as required.\qed \\

If $S$ is an inverse semigroup then $\mathsf{P}(S)$ is a semigroup with involution when we define
$AB = \{ab \colon a \in A, b \in B \}$ and $A^{-1} = \{a^{-1} \colon a \in A \}$.
The structure of filters on boolean inverse monoids is closely bound up with the following definition.
A {\em coset} in $S$ is a subset $A$ such that $A = AA^{-1}A$.
The theory of upwardly closed cosets in inverse semigroups is discussed in detail in \cite{L1}.

\begin{lemma} Every filter is a coset.
\end{lemma}
\proof Let $F$ be a filter and $ab^{-1}c \in FF^{-1}F$ where $a,b,c \in F$.
Put $d = a \wedge b \wedge c$.
Then $d = dd^{-1}d \leq ab^{-1}c$ and so $ab^{-1}c \in F$.
The reverse inclusion is immediate.\qed\\

\begin{lemma} If $A$ and $B$ are filters then $(AB)^{\uparrow}$ is the smallest filter containing $AB$.
\end{lemma}
\proof Clearly $AB \subseteq (AB)^{\uparrow}$ and $(AB)^{\uparrow}$ is upwardly closed.
We show that $(AB)^{\uparrow}$ is a filter base.
Let $a_{1},a_{2} \in A$ and $b_{1},b_{2} \in B$.
Then $a = a_{1} \wedge a_{2} \in A$ and $b = b_{1} \wedge b_{2} \in B$.
Now $ab \leq a_{1}b_{1}$ and $ab \leq a_{2}b_{2}$.
Thus $ab \leq a_{1}b_{1} \wedge a_{2}b_{2}$ and $ab \in AB$.
It follows that $a_{1}b_{1} \wedge a_{2}b_{2} \in (AB)^{\uparrow}$.\qed\\

In the light of the above lemma, we may define the filter $A \cdot B = (AB)^{\uparrow}$ when $A$ and $B$ are filters.
The following is Proposition~1.4 of \cite{L1}.

\begin{lemma} Let $F$ be a filter. Then $H = F^{-1} \cdot F$ is a filter and an inverse submonoid and $F = (aH)^{\uparrow}$ for any $a \in F$.
\end{lemma}

An {\em idempotent filter} is a filter containing an idempotent.
The following is Proposition~1.5 of \cite{L1}.

\begin{lemma} A filter $F$ is idempotent if and only if it is an inverse subsemigroup.
\end{lemma}

\begin{remark}{\em  Our results on filters are special cases of some well-known results on actions of inverse semigroups
and their associated closed inverse subsemigroups ultimately due to Boris Schein.
If $H$ is an idempotent filter then the {\em left cosets} of $H$ are the sets of the form $(aH)^{\uparrow}$ where $\dom (a) \in H$;
it is important to take note of the extra condition which is, of course, automatic in the case of groups since they only have one idempotent.
Observe that $(aH)^{\uparrow} = (bH)^{\uparrow}$ iff $a^{-1}b \in H$.}
\end{remark}

The following result is a corollary to Lemma~2.8.

\begin{lemma} Let $A$ and $B$ be filters such that $A \cap B \neq \emptyset$ and $A^{-1} \cdot A = B^{-1} \cdot B$ then $A = B$.
\end{lemma}

Denote by $\mathsf{L}(S)$ the set of all filters on $S$ equipped with the product $\cdot$ defined above.
Either directly or via \cite{Lenz} we have the following.

\begin{proposition} Let $S$ be a boolean inverse monoid. 
Then $\mathsf{L}(S)$ is an inverse semigroup in which the idempotents are the idempotent filters and the natural partial order is reverse inclusion.
\end{proposition}

Denote by $\mathsf{G}(S)$ the subset of $\mathsf{L}(S)$ consisting of all the ultrafilters of $S$.

\begin{proposition} Let $S$ be a boolean inverse monoid. 
Then with respect to the semigroup multiplication in  $\mathsf{L}(S)$, the set  $\mathsf{G}(S)$ is a groupoid. 
In addition, the following are equivalent
\begin{enumerate}

\item $F$ is an ultrafilter in $S$.

\item $H = F^{-1} \cdot F$ is an idempotent ultrafilter in $S$.

\item  $E(H)$ is an ultrafilter in $E(S)$.
\end{enumerate}

\end{proposition}
\proof An element $s \in S$ in an inverse monoid $S$ is said to be {\em primitive} if $t \leq s$ and $t \neq s$ implies that $t = 0$.
The primitive elements of an inverse monoid form a groupoid by Proposition~9.2.1 of \cite{L}.
Because the natural partial order in $\mathsf{L}(S)$ is reverse inclusion,
it follows that the primitive elements in $\mathsf{L}(S)$ are precisely the ultrafilters in $S$.

The equivalence of (1) and (2) follows from the above, but we prove it directly anyway.

(1)$\Rightarrow$(2). Put $H = F^{-1} \cdot F$.
Let $H \subseteq K$ where $K$ is a filter.
Since $H$ is an idempotent filter so too is $K$ by Lemma~2.9.
Let $a \in F$.
Then $\dom (a) \in H$ and so $\dom (a) \in K$.
Thus $F = (aH)^{\uparrow} \subseteq (aK)^{\uparrow}$.
But by assumption $F$ is an ultrafilter and so $F = (aH)^{\uparrow} = (aK)^{\uparrow}$.
It follows that $H = K$ and so $H$ is also an ultrafilter.

(2)$\Rightarrow$(3). Put $E' = E(F^{-1} \cdot F)$.
Let $E' \subseteq F'$ where $F'$ is a filter in $E(S)$.
Then $H \subseteq F'^{\uparrow}$. 
But $H$ is an ultrafilter and so $H = F'^{\uparrow}$. 
Thus $E' = F'$, and so $E'$ is an ultrafilter in the semilattice of idempotents.
 
(3)$\Rightarrow$(1). We have that $F = (aH)^{\uparrow}$ where $E(H)$ is an ultrafilter in $E(S)$.
Suppose that $F \subseteq G$ where $G$ is a filter.
Then $H \subseteq G^{-1} \cdot G$ and so $E(H) \subseteq E(G^{-1} \cdot G)$.
By assumption  $E(H) = E(G^{-1} \cdot G)$ and so $H = G^{-1} \cdot G$ 
from which it follows that $F = G$ and so $F$ is an ultrafilter.\qed \\

We can easily write down an explicit form of the groupoid multiplication in $\mathsf{G}(S)$.
Let $A$ and $B$ be two ultrafilters such that $A^{-1} \cdot A = B \cdot B^{-1}$
Then $A = (aA^{-1} \cdot A)^{\uparrow}$ where $a \in A$ and $B = (B \cdot B^{-1} b)^{\uparrow}$ where $b \in B$.
Thus $A \cdot B = (ab B^{-1} \cdot B)^{\uparrow}$ since $ab \in A \cdot B$ and 
$(A \cdot B)^{-1} \cdot (A \cdot B) = B^{-1} \cdot B$.

\begin{lemma} Let $A$ be an ultrafilter in a boolean inverse monoid $S$.
Then $s \vee t \in A$ implies that $s \in A$ or $t \in A$.
\end{lemma}
\proof Put $H = A^{-1} \cdot A$. Then $A = (aH)^{\uparrow}$ for any $a \in A$.
Now $s \vee t \in A$ implies that $\dom (s \vee t) \in H$.
Thus $\dom (s) \vee \dom (t) \in E(H)$.
However $E(H)$ is an ultrafilter in the boolean algebra $E(S)$.
Thus $\dom (s) \in H$ or $\dom (t) \in H$.
Suppose the former.
Put $B = (sH)^{\uparrow}$, a well-defined ultrafilter.
Observe that $s \in B$ implies that $s \vee t \in B$.
Then $B^{-1} \cdot B = A^{-1} \cdot A$ and $s \vee t \in A \cap B \neq \emptyset$.
By Lemma~2.11 it follows that $A = B$ and so $s \in A$, as required.\qed \\

\begin{proposition} There is a contravariant functor $\mathsf{G}$ from the category of boolean inverse monoids 
to the category of groupoids and their covering functors.
\end{proposition}
\proof Let $\theta \colon S \rightarrow T$ be a morphism between boolean inverse monoids.
By (M3), the function $\theta^{-1} \colon \mathsf{G}(T) \rightarrow \mathsf{G}(S)$ is well-defined.
It is easy to check that $\theta^{-1}(A)^{-1} = \theta (A^{-1})$.

We prove first that this function is a functor.
Let $x \in (\theta^{-1}(A^{-1}) \theta^{-1}(A))^{\uparrow}$. 
Then $uv \leq x$ where $u \in \theta^{-1}(A^{-1})$ and $v \in \theta^{-1}(A))^{\uparrow}$. 
Thus $\theta (u) \in A^{-1}$ and $\theta (v) \in A$ and so $\theta (uv) \in A^{-1}A$.
But $\theta (uv) \leq \theta (x)$.
Thus $\theta (x) \in (A^{-1}A)^{\uparrow}$.
We have therefore proved that
$$(\theta^{-1}(A^{-1}) \theta^{-1}(A))^{\uparrow} \subseteq \theta^{-1}((A^{-1}A)^{\uparrow}).$$
But both sets are ultrafilters and so this inclusion must be equality.
We have therefore proved that 
$$\theta^{-1}(A^{-1} \cdot A) = \theta^{-1}(A)^{-1} \cdot \theta^{-1}(A).$$
The dual result holds by symmetry.
It follows that $\theta^{-1}$ preserves domains and ranges.
To conclude the proof that $\theta^{-1}$ is a functor we prove two results and then combine then.

First, we prove that $\theta^{-1}(A) \theta^{-1}(B) \subseteq \theta^{-1}(AB)$.
Let $x \in \theta^{-1}(A)\theta^{-1}(B)$.
Then $x = a'b'$ where $a' \in \theta^{-1}(A)$, $b' \in \theta^{-1}(B)$.
Thus $\theta (a') \in A$ and $\theta (b') \in B$ giving $\theta (x) = \theta (a'b')$.
Hence $\theta (x) \in \theta^{-1}(AB)$. 

Second, we prove that $\theta^{-1} (X)^{\uparrow} \subseteq \theta^{-1}(X^{\uparrow})$.
Let $a \in \theta^{-1}(X)^{\uparrow}$.
Then $b \leq a$ where $b \in \theta^{-1}(X)$.
Thus $\theta (b) \in X$ and so $\theta (b) \leq \theta (a)$.
Thus $\theta (a) \in X^{\uparrow}$ giving $a \in \theta^{-1}(X^{\uparrow})$.

We now combine these two results.
We have that 
$$\theta^{-1}(A) \theta^{-1}(B) \subseteq \theta^{-1}(AB).$$
Thus 
$$(\theta^{-1}(A) \theta^{-1}(B))^{\uparrow} 
\subseteq \theta^{-1}(AB)^{\uparrow} 
\subseteq \theta^{-1} ( (AB)^{\uparrow} ).$$
This gives $\theta^{-1}(A) \cdot \theta^{-1}(B) \subseteq \theta^{-1}(A \cdot B)$.
But again we have ultrafilters on both sides of the inclusion and so it is in fact an equality.
This proves that $\theta^{-1}$ is a functor.

We prove that $\theta^{-1}$ is a covering functor.
Let $A,B \in \mathsf{G}(T)$ such that $A^{-1} \cdot A = B^{-1} \cdot B$ and $\theta^{-1}(A) = \theta^{-1} (B)$.
Then $A \cap B \neq \emptyset$.
Thus $A = B$ by Lemma~2.11.
It follows that $\theta^{-1}$ is star injective.

Let $H \in \mathsf{G}(T)$ be an idempotent ultrafilter such that $\theta^{-1}(H) = B^{-1} \cdot B$ where $B$ is an ultrafilter in $S$.
Let $b \in B$ and put $a = \theta (b)$.
Then $a^{-1}a = \theta (b^{-1}b)$.
Now $b^{-1}b \in B^{-1} \cdot B$ and so $b^{-1}b \in \theta^{-1}(H)$ and so $\theta (b^{-1}b) \in H$.
Thus $a^{-1}a \in H$.
It follows that $(aH)^{\uparrow}$ is a well-defined ultrafilter in $T$ whose domain is $H$.
Observe that $\theta^{-1} (aH) \cap B \neq \emptyset$.
Thus $\theta^{-1}(aH) = B$ by Lemma~2.11.
It follows that $\theta^{-1}$ is star surjective.\qed \\

Let $G$ be a groupoid.
The set of bisections of $G$, denoted by $\mathsf{B}(G)$, forms an inverse semigroup \cite{Renault,Paterson};
in fact, it is a boolean inverse monoid.
The following is a well-known property of covering functors between groupoids.

\begin{lemma} Let $\alpha \colon G \rightarrow H$ be a covering functor between two groupoids.
Then if $\theta (x) = ab$ there exist $u,v \in G$ such that $x = uv$,
$\alpha (u) = a$ and $\alpha (v) = b$. 
\end{lemma}

The process of constructing the inverse monoid of bisections of a groupoid is functorial.

\begin{proposition} There is a contravariant functor $\mathsf{B}$ from the category of groupoids and their covering functors 
to the category of boolean inverse monoids and their monoid homomorphisms that preserve meets.
\end{proposition}
\proof Let $\alpha \colon G \rightarrow H$ be a covering functor between groupoids.
Let $B$ be a bisection of $H$.
Let $a,b \in \alpha^{-1}(B)$ such that $\dom (a) = \dom (b)$.
Then $\alpha (a), \alpha (b) \in B$ and $\dom (\alpha (a)) = \dom (\alpha (b))$.
But $B$ is a bisection and so $\alpha (a) = \alpha (b)$.
But $\alpha $ is star-injective and so $a = b$.
Together with a dual argument, this proves that $\alpha^{-1}(B)$ is a bisection.
We therefore have a well-defined function $\alpha^{-1} \colon \mathsf{B}(H) \rightarrow \mathsf{B}(G)$.
This map induces a homomorphism between the boolean algebras of idempotents and $\alpha^{-1}(A \cap B) = \alpha^{-1}(A) \cap \alpha^{-1}(B)$.
We prove that $\alpha^{-1}(AB) = \alpha^{-1}(A)\alpha^{-1}(B)$ and so $\alpha^{-1}$ is a homomorphism.
Let $x \in \alpha^{-1}(AB)$.
Then $\alpha (x) = ab$.
Then by Lemma~2.16 there exists $u,v$ such that $x = uv$ and $\alpha (u) = a$ and $\alpha (v) = b$.
Thus $u \in \alpha^{-1}(A)$ and $v \in \alpha^{-1}(B)$ and so $x \in \alpha^{-1}(A)\alpha^{-1}(B)$.
Thus we have proved that  $\alpha^{-1}(AB) \subseteq \alpha^{-1}(A)\alpha^{-1}(B)$.
To prove the reverse inclusion let $x \in  \alpha^{-1}(A)\alpha^{-1}(B)$.
Then $x = uv$ where $\alpha (u) \in A$ and $\alpha (v) \in B$.
Thus $\alpha (x) \in AB$ and so $x \in \alpha^{-1}(AB)$, as required. 
We may therefore define $\mathsf{B}(\alpha) = \alpha^{-1}$.\qed \\

So far we have only dealt with matters algebraical, we now deal with those topological.
Let $G$ be a boolean groupoid.
Denote by $\mathsf{A}(G)$ the set of compact open bisections of $G$.
We could take some shortcuts in the proof of the proposition below using \cite{Paterson}
but we have tried to be as elementary and explicit as possible.
Observe by (6) below, that boolean groupoids are {\em locally compact}.

\begin{proposition} Let $G$ be a boolean groupoid.
\begin{enumerate}

\item $G_{o}$ is an open set in $G$.

\item  $G \ast G$ is a closed set in $G \times G$.

\item If $A$ is a closed bisection then $A^{-1}A$ is a closed subset of $G_{o}$.

\item The product of two open sets is an open set.

\item The product of two closed sets is a closed set.

\item An open bisection is closed if and only if it is compact.
Thus clopen bisections are the same thing as compact open bisections.

\item The product of two compact open bisections is a compact open bisection.

\item  $\mathsf{A}(G)$ is a boolean inverse monoid.

\end{enumerate}
\end{proposition}
\proof (1) In an \'etale topological groupoid $G_{o}$ is always open \cite{R1}.

(2) In a topological groupoid $G$ in which $G_{0}$ is hausdorff the set $G \ast G$ is closed.
This follows from general topology.

(3) Because $A$ is a bisection, $A^{-1}A = \{\dom(a) \colon a \in A \} = \dom (A)$, a subset of $G_{o}$.
In an \'etale topological groupoid, the map $\dom  \colon G \rightarrow G_{o}$ is open \cite{R1}, and so it maps closed sets to closed sets.
Thus if $A$ is a closed bisection, then $A^{-1}A$ is a closed subset of $G_{o}$.

(4) In an \'etale topological groupoid the product of any two open sets is an open set \cite{R1}.

(5) Let $A$ and $B$ be closed sets in $G$.
Then $A \times B$ is a closed set in $G \times G$.
Because our topological groupoid is hausdorff, (2) implies that $A \ast B = (A \times B) \cap (G \ast G)$ is a closed subset of $G \times G$.
In an \'etale topological groupoid the multiplication map is open and so maps closed sets to closed sets.
Thus $AB$ is closed.

(6) Observe that to determine whether a subset of a space is compact it is enough to use covers whose elements are taken from a basis for the topology.
In our case, the basis consists of compact open bisections.
Let $A$ be a clopen bisection which we are required to show is compact.
Let $A = \bigcup B_{i}$ be an open cover of $A$, where by our observation above each $B_{i}$ is a compact open bisection.
Now the $B_{i}$ are pairwise compatible since they are bounded above,
and so by Lemma~2.1(2), $A^{-1}A = \bigcup B_{i}^{-1}B_{i}$ where $A^{-1}A$ and the $B_{i}^{-1}B_{i}$ are subsets of $G_{o}$.
But by (2) above, $A^{-1}A$ is a closed subspace of a compact space $G_{o}$ and so it is compact.
It follows that we can write $A^{-1}A = \bigcup_{i=1}^{m} B_{i}^{-1}B_{i}$ for some finite number of elements.
Now $\bigcup_{i=1}^{m} B_{i} \subseteq A$ and the domains of $\bigcup_{i=1}^{m} B_{i}$ and $A$ are the same. 
It follows that $A = \bigcup_{i=1}^{m} B_{i}$, and so we have proved that $A$ is compact.
Conversely, let $A$ be a compact open bisection.
But every compact subset of a hausdorff space is closed and so $A$ is clopen.

(7) This follows by (4),(5) and (6) and the fact that the product of bisections is a bisection.

(8) The inverse of a compact open bisection is a compact open bisection and so together with (6)
we have that $\mathsf{A}(G)$ is an inverse subsemigroup of $\mathsf{B}(G)$.
By (1), $G_{0}$ is an open subspace of $G$ \cite{R1}, it is compact by {\em fiat}, and it is automtaically a bisection. 
Thus $G_{o}$ is a compact open set and is the identity element of $\mathsf{A}(S)$.
Now $G$ has a basis of compact open bisections each of which is clopen, 
thus these intersect with $G_{0}$ to give a basis of clopen subsets for $G_{o}$.
It follows that $G_{0}$ is a boolean space and this implies that the idempotents of $\mathsf{A}(G)$ form a boolean algebra.
The natural partial order in $\mathsf{A}(G)$ is just subset inclusion.
If $A,B \in \mathsf{A}(G)$ then $A \cap B$ is a clopen bisection.
It follows that $\mathsf{A}(G)$ has all non-empty finite meets. 
Finally, if $A,B \in \mathsf{A}(G)$ are orthogonal then $A \cup B$ is a bisection and it is clopen since both $A$ and $B$ are clopen.
Hence $A \cup B \in \mathsf{A}(G)$ and so $\mathsf{A}(G)$ has finite orthogonal joins. 
Thus $\mathsf{A}(G)$ is a boolean inverse monoid, as claimed.\qed\\

Let $G$ be a boolean groupoid.
For each $g \in G$, 
define 
$$\mathcal{F}_{g} = \{A \in \mathsf{A}(G) \colon g \in A \}.$$

\begin{lemma} With the above definition, we have the following.
\begin{enumerate}

\item $\mathcal{F}_{g}$ is an ultrafilter in the inverse semigroup $\mathsf{A}(G)$.

\item $\mathcal{F}_{g}^{-1} \cdot \mathcal{F}_{g} = \mathcal{F}_{g^{-1}g}$ and  
$\mathcal{F}_{g} \cdot \mathcal{F}_{g}^{-1} = \mathcal{F}_{gg^{-1}}$ 

\item If $\exists gh$ then $\mathcal{F}_{g} \cdot \mathcal{F}_{h} = \mathcal{F}_{gh}$.

\item $\mathcal{F}_{g} = \mathcal{F}_{h}$ iff $g = h$.

\item Each ultrafilter $F$ in the boolean inverse monoid $\mathsf{A}(G)$ is of the form $\mathcal{F}_{g}$ for some $g \in G$.

\end{enumerate}
\end{lemma}
\proof (1) It is immediate that $\mathcal{F}_{g}$ is a filter. It remains to show that it is an ultrafilter.
Let $A \in \mathsf{A}(G)$ be a compact open bisection with the property that for each $B \in \mathcal{F}_{g}$ we have that $A \cap B \neq \emptyset$.
We shall show that $g \in A$ from which it follows that $A \in \mathcal{F}_{g}$ and so by Lemma~2.4 we deduce that $\mathcal{F}_{g}$ is an ultrafilter.
Let $O$ be any open set containing $g$.
Then there is a basic compact open set $g \in D$ such that $D \subseteq O$.
By assumption $D \cap A \neq \emptyset$ and so $O \cap A \neq \emptyset$.
Thus every open set containing $g$ intersects $A$ non-emptily.
But $A$ is closed and so $g \in A$, as required.

(2) It is clear that $\mathcal{F}_{g}^{-1} \cdot \mathcal{F}_{g} \subseteq \mathcal{F}_{g^{-1}g}$.
But we now use the fact that the lefthand side is an ultrafilter and the righthand side an (ultra)filter.
They must therefore be equal.

(3) Similar argument to (2) above.

(4) Suppose that $g \neq h$.
By assumption the groupoid $G$ is hausdorff.
Thus there are basic compact open bisections $A$ and $B$ such that $g \in A$ and $h \in B$ and $A \cap B = \emptyset$
which shows that $\mathcal{F}_{g} \neq \mathcal{F}_{h}$.

(5) Let $F$ be an ultrafilter in $\mathsf{A}(G)$ 
We prove that $F \subseteq \mathcal{F}_{g}$ for some $g \in G$ from which the result follows.
Let $A \in F$.
Then $A$ is a compact set.
Consider the set $F' = \{ A \cap B \colon B \in F \}$.
Then this is a set of closed subsets of $A$ which has the finite intersection property because $F$ is an ultrafilter.
It follows that $\bigcap F'$ is non-empty.
Let $g \in \bigcap F'$.
Then $g$ belongs to every element of $F$
and so by construction $F \subseteq \mathcal{F}_{g}$, as required.\qed \\

\begin{proposition} The construction $\mathsf{A}$ is a contravariant functor from the category of boolean groupoids 
to the category of boolean inverse monoids.
\end{proposition}
\proof By Proposition~2.18, $\mathsf{A}(G)$ is a boolean inverse monoid.
Let $\alpha \colon G \rightarrow H$ be a continuous proper covering functor.
Thus $\alpha^{-1}$ takes clopen bisections of $H$ to clopen bisections of $G$
and so compact open bisections to compact open bisections.
Combining this observation with Proposition~2.17, 
we have that $\alpha^{-1} \colon \mathsf{A}(H) \rightarrow \mathsf{A}(G)$ is a monoid homomorphism which preserves meets.
It is clearly a boolean algebra map from the boolean algebra of idempotents of $\mathsf{A}(H)$ to the boolean algebra of idempotents of $\mathsf{A}(G)$.
It remains to prove that the inverse images of ultrafilters in $\mathsf{A}(G)$ are ultrafilters in $\mathsf{A}(H)$.
Let $F$ be an ultrafilter in $\mathsf{A}(G)$.
Then by Lemma~2.19(5), there exists $g \in G$ such that $F = \mathcal{F}_{g}$.
Put $h = \alpha (g) \in H$.
Then $\mathcal{F}_{h}$ is an ultrafilter in $\mathsf{A}(H)$.
The result will be proved if we can show that
$A \in \mathcal{F}_{h} \Leftrightarrow \alpha^{-1}(A) \in \mathcal{F}_{g}$.
Suppose that $A \in  \mathcal{F}_{h}$.
Then $h \in A$ and so $\alpha (g) \in A$ which gives that $g \in \alpha^{-1}(A)$ and so $\alpha^{-1}(A) \in \mathcal{F}_{g}$.
Conversely, suppose that $\alpha^{-1}(A) \in \mathcal{F}_{g}$.
Then $g \in \alpha^{-1}(A)$ and so $h \in A$ giving $A \in \mathcal{F}_{h}$.\qed \\

Let $S$ be a boolean inverse monoid.
For each $s \in S$ define
$$\mathcal{K}_{s} = \{A \in \mathsf{G}(S) \colon s \in A \}.$$
It follows by (1) below that the set $\{\mathcal{K}_{s} \colon s \in S \}$ is a basis for a topology $\Omega$ on the groupoid $\mathsf{G}(S)$.

\begin{lemma} Let $S$ be a boolean monoid. 
With the above definition we have the following.
\begin{enumerate}

\item Each $\mathcal{K}_{s}$ is a bisection.

\item $\mathcal{K}_{s} \cap \mathcal{K}_{t} = \mathcal{K}_{s \wedge t}$.

\item $\mathcal{K}_{s}^{-1} = \mathcal{K}_{s^{-1}}$.

\item $\mathcal{K}_{s}\mathcal{K}_{t} = \mathcal{K}_{st}$.

\item $\mathcal{K}_{s} \subseteq \mathcal{K}_{t}$ iff $s \leq t$.

\item  $\mathcal{K}_{s} = \mathcal{K}_{t}$ iff $s = t$. 

\item If $s \vee t$ exists then $\mathcal{K}_{s} \cup \mathcal{K}_{t} = \mathcal{K}_{s \vee t}$.

\item $\mathcal{K}_{s} \cup \mathcal{K}_{t}$ is a bisection iff $s \vee t$ exists.

\item Each compact open bisection of $\mathsf{G}(S)$ is equal to $\mathcal{K}_{s}$ for some $s \in S$.

\end{enumerate}
\end{lemma}
\proof (1) Let $A,B \in \mathcal{K}_{s}$ such that $A^{-1} \cdot A = B^{-1} \cdot B$.
Then $s \in A \cap B$ and so in particular $A \cap B \neq \emptyset$.
Thus by Lemma~2.11, we have that $A = B$.
The dual result also holds. 

(2) Straightforward.

(3) Let $A \in \mathcal{K}_{s} \cap \mathcal{K}_{t}$.
Then $s,t \in A$.
But $A$ is a filter and so $s \wedge t \in A$ from which it follows that $A \in \mathcal{K}_{s \wedge t}$.
Conversely, let $A \in \mathcal{K}_{s \wedge t}$.
Then $s \wedge t \in A$.
But $s \wedge t \leq s,t$ and $A$ is a filter and so $s,t \in A$ and so $A \in \mathcal{K}_{s} \cap \mathcal{K}_{t}$.
 
(4) Let $A \in \mathcal{K}_{s}$ and $B \in \mathcal{K}_{t}$.
Then $st \in A \cdot B$ and so $\mathcal{K}_{s} \cap \mathcal{K}_{t} \subseteq \mathcal{K}_{s \wedge t}$.
Conversely, let $A \in \mathcal{K}_{st}$.
Put $H = A^{-1} \cdot A$.
Then $A = (stH)^{\uparrow}$.
Put $B = (s(tHt^{-1})^{\uparrow})^{\uparrow}$ and $C = (tH)^{\uparrow}$.
Then $B \in \mathcal{K}_{s}$, and $C \in \mathcal{K}_{t}$ and $A = B \cdot C$.
Thus $\mathcal{K}_{s \wedge t} \subseteq \mathcal{K}_{s} \cap \mathcal{K}_{t}$.
 
(5) Let $\mathcal{K}_{s} \subseteq \mathcal{K}_{t}$.
Suppose that $s \nleq t$.
Then by Lemma~2.2(3), there exists a non-zero element $s'$ such that $s' \leq s$ and $s' \wedge t = 0$.
By Lemma~2.5(1), let $A$ be an ultrafilter containing $s'$.
Then $A$ contains $s$.
Thus $A \in \mathcal{K}_{s}$ and so $A \in \mathcal{K}_{t}$.
It follows that $t \in A$.
But $s',t \in A$ implies that $s' \wedge t \in A$ and so $0 \in A$, which is a contradiction.
It follows that $s \leq t$.
The converse is immediate since any ultrafilter containing $s$ must contain $t$ as well.

(6) This is immediate by (5).

(7) This follows by Lemma~2.14.

(8) One direction is immediate.
Suppose that $\mathcal{K}_{s} \cup \mathcal{K}_{t}$ is a bisection.
We prove that $s$ and $t$ are compatible.
We prove that $st^{-1}$ is an idempotent.
The fact that $s^{-1}t$ is an idempotent follows by symmetry.
If $st^{-1} = 0$ there is nothing to prove so we may assume that $st^{-1} \neq 0$.
Thus $e = s^{-1}st^{-1}t \neq 0$.
Let $H$ be any ultrafilter containing $e$; such exists by Lemma~2.5(1).
Necessarily $H$ is an idempotent ultrafilter by Lemma~2.9.
Since $e \in H$ both $s^{-1}s,t^{-1}t \in H$.
Put $A = (sH)^{\uparrow}$ and $B = (tH)^{\uparrow}$.
Then $A \in \mathcal{K}_{s}$ and $B \in \mathcal{K}_{t}$.
But $A^{-1} \cdot A = B^{-1} \cdot B$ by construction.
It follows that $A = B$.
Thus $s^{-1}t \in H$.
We have proved that every ultrafilter containing $e$ contains $s^{-1}t$.
By Lemma~2.5(2), we have that $s^{-1}st^{-1}t \leq s^{-1}t$.
Thus $st^{-1} \leq ss^{-1}tt^{-1}$.
We have proved that $st^{-1}$ is an idempotent.

(9) Let $A$ be a compact open bisection in the groupoid $\mathsf{G}(S)$.
Because it it an open set it is a union of basic compact open bisections,
and because it is compact it is a union of only a finite number of these sets.
Thus 
$$A = \bigcup_{i=1}^{m} \mathcal{K}_{s_{i}}.$$
By (8) above the elements $s_{1}, \ldots, s_{m}$ are pairwise compatible.
Thus by Lemma~2.3, the join $s = \bigvee_{i=1}^{m} s_{i}$ exists.
It follows by (7) above that $A = \mathcal{K}_{s}$, as required.\qed \\

\begin{proposition} Let $S$ be a boolean inverse monoid.
Then $\mathsf{G}(S)$ is a boolean groupoid with respect to the topology $\Omega$.
If $\theta \colon S \rightarrow T$ is a morphism of boolean inverse monoids then
$\mathsf{G}(\theta) \colon \mathsf{G}(S) \rightarrow \mathsf{G}(T)$ is
a continuous covering functor.
Thus the construction $\mathsf{G}$ is a contravariant functor from the category of boolean inverse monoids to the category of boolean groupoids.
\end{proposition}
\proof We procede in a number of steps.

The basic open sets are also closed.
Let $s \in S$ be a non-zero element and let $F \in \mathsf{G}(S) \setminus \mathcal{K}_{s}$.
By assumption $s \notin F$ and so by Lemma~2.4, there exists $t \in F$ such that $s \wedge t = 0$.
Thus $F \in \mathcal{K}_{t}$ and $\mathcal{K}_{t} \cap \mathcal{K}_{s} = \emptyset$.
It follows that $\mathcal{K}_{s}$ is also a closed subset.

The topology $\Omega$ is hausdorff.
Let $F$ and $G$ be two distinct ultrafilters in $S$; in other words elements of the groupoid $\mathsf{G}(S)$.
If $s \in F$ and $s \wedge G \neq 0$ then $s \in G$ by Lemma~2.4.
Since $F$ cannot be a subset of $G$ there must exist $s \in F$ and $t \in G$ such that $s \wedge t = 0$.
Then $F \in \mathcal{K}_{s}$, $G \in \mathcal{K}_{t}$ and $\mathcal{K}_{s} \cap \mathcal{K}_{t} = \emptyset$.

$\mathsf{G}(S)$ is a topological groupoid.
We have to prove that the inversion map and the multiplication map are both continuous.
The fact the the inversion map is continuous follows by Lema~2.21(3).
We prove that the multiplication map $\mu \colon \mathsf{G}(S) \ast \mathsf{G}(S) \rightarrow \mathsf{G}$ is continuous.
To do this we prove that
$$\mu^{-1}(\mathcal{K}_{a}) = (\bigcup_{0 \neq bc \leq a} \mathcal{K}_{b} \times \mathcal{K}_{c}) \cap (\mathsf{G}(S) \ast \mathsf{G}(S).$$
Let $a \in A$ be an ultrafilter such that $A = B \cdot C$.
Then $a \in (BC)^{\uparrow}$ and so $bc \leq a$ for some $b \in B$ and $c \in C$.
Thus $B \in \mathcal{K}_{b}$, $C \in \mathcal{K}_{c}$ and $0 \neq bc \leq a$.
To prove the reverse inclusion, 
suppose that $0 \neq bc \leq a$ and $B \in \mathcal{K}_{b}$, $C \in \mathcal{K}_{c}$ and the product $B \cdot C$ exists.
Then $B \cdot C$ is an ultrafilter containing $a$ and so $B \cdot C \in \mathcal{K}_{a}$.

$\mathsf{G}(S)$ is \'etale. We shall show that $\dom \colon \mathsf{G}(S) \rightarrow \mathsf{G}(S)_{o}$ is a local homeomorphism.
To do this it is enough to prove that the map $\dom \colon \mathcal{K}_{s} \rightarrow \mathcal{K}_{s^{-1}s}$
given by $A \mapsto A^{-1} \cdot A$ is a homeomorphism.
It is bijective by Lemmas~2.8 and 2.11.
It is continuous because inversion and multiplication are continuous.
To show that it is open, we use Lemma~2.21(5):
$\mathcal{K}_{t}$ is an open set in $\mathcal{K}_{s}$ iff $t \leq s$.
It follows that $\mathcal{K}_{t^{-1}t}$ is an open set in $\mathcal{K}_{s^{-1}s}$.

The fact that $\mathsf{G}(S)_{o}$ is compact follows from the proof of classical Stone duality.
We have therefore proved that $\mathsf{G}(S)$ is a boolean groupoid. 

By Proposition~2.15, it only remains to show that $\theta^{-1}$ is continuous.
Let $\mathcal{K}_{s}$ be a basic open set in  $\mathsf{G}(S)$.
Put $t = \theta (s)$.
Then $\mathcal{K}_{t}$ is a basic open set in $\mathsf{G}(T)$.
Our claim will be proved if we can show that
$F \in \mathcal{K}_{t} \Leftrightarrow \theta^{-1}(F) \in \mathcal{K}_{s}$.
But the proof of this is immediate.\qed \\

The following concludes the proof of our main result by showing that our contravariant functors $\mathsf{A}$ and $\mathsf{G}$ 
establish a dual equivalence between the categories of boolean inverse monoids and boolean groupoids.

\begin{proposition} \mbox{}

\begin{enumerate}

\item Let $G$ be a boolean groupoid.
Then $G$ is isomorphic to $\mathsf{G}\mathsf{A}(G)$ as topological groupoids under the map $g \mapsto \mathcal{F}_{g}$.

\item Let $S$ be a boolean inverse monoid.
Then $S$ is isomorphic to  $\mathsf{A}\mathsf{G}(S)$ under the map $s \mapsto \mathcal{K}_{s}$.

\end{enumerate}
\end{proposition}
\proof (1) By Lemma~2.19(4) and (5), the groupoids are isomorphic under this map so it only remains to prove that this map is a homeomorphism.
We show first that the map is open.
Let $U$ be a compact open bisection of $G$.
We show that
$$\{\mathcal{F}_{g} \colon g \in U \} = \mathcal{K}_{U}.$$
Let $g \in U$.
Then $\mathcal{F}_{g}$ is an ultrafilter of $\mathsf{A}(G)$.
It consists of all compact open bisections of $G$ that contain $g$.
But $U$ is such a one.
Thus $U \in \mathcal{F}_{g}$ and so $\mathcal{F}_{g} \in \mathcal{K}_{U}$.
Now let $A \in \mathcal{K}_{U}$.
Then $A$ is an ultrafilter of $\mathsf{A}(G)$ and $U \in A$.
But by Lemma~2.19(5), we know that $A = \mathcal{F}_{h}$ for some $h \in G$.
But $U \in A$ and so $h \in U$.
It follows that $A \in \{\mathcal{F}_{g} \colon g \in U \}$.
Finally we prove that it is continuous.
A basic compact open subset of $\mathsf{G}\mathsf{A}(G)$ is of the form $\mathcal{K}_{U}$ where $U$ is a compact open bisection in $G$.
The inverse image of $\mathcal{K}_{U}$ under the map is the set of all elements $g \in G$ such that $\mathcal{F}_{g} \in \mathcal{K}_{U}$.
But $\mathcal{F}_{g} \in \mathcal{K}_{U}$ iff $U \in \mathcal{F}_{g}$ iff $g \in U$.
Thus the inverse image of $\mathcal{K}_{U}$ under the map is $U$.

(2) This is immediate by Lemma~2.21(4) and (9).\qed \\

\begin{remark}{\em The one aspect of our duality that is not as straightforward as we would have hoped is 
our definition of a morphism $\theta \colon S \rightarrow T$ between boolean inverse monoids.
This includes the requirement that the inverse images of ultrafilters be ultrafilters.
We can see why some such condition is needed by considering the case where $S$ is a boolean algebra.
Then the inverse image of any non-idempotent ultrafilter in $T$ has to be empty.
If we waive (M3) and assume only (M1) and (M2) we can prove that the inverse images of idempotent ultrafilters
are idempotent using Lemma~2.13 as follows.
Let $H \subseteq T$ is an idempotent ultrafilter.
Then $H \cap E(T)$ is an ultrafilter in the boolean algebra $E(T)$.
Thus $\theta^{-1}(H \cap E(T)) \cap E(S)$ is an ultrafilter in the boolean algebra $E(S)$.
It follows that $(\theta^{-1}(H \cap E(T)) \cap E(S))^{\uparrow}$ is an idempotent ultrafilter in $S$ 
But  $(\theta^{-1}(H \cap E(T)) \cap E(S))^{\uparrow} \subseteq \theta^{-1}(H)$ and $\theta^{-1}(H)$ is a filter.
so that $(\theta^{-1}(H \cap E(T)) \cap E(S))^{\uparrow} = \theta^{-1}(H)$
showing that $\theta^{-1}(H)$ is an idempotent ultrafilter.}
\end{remark}

We conclude this section with three examples.

\begin{example}{\em Our first example shows that our theory is a proper extension of Stone duality.
Let $X$ be a finite non-empty set.
Then the symmetric inverse monoid, $I(X)$, on $X$ is a boolean inverse monoid;
the semilattice of idempotents of $I(X)$ is isomorphic to the boolean algebra of all subsets of $X$.
Because $I(X)$ is finite and has finite intersections, each filter in $I(X)$ is principal
and the ultrafilters are in bijective correspondence with the elements of $I(X)$ 
whose domains, and therefore whose ranges, contain exactly one element.
The boolean groupoid associated with $I(X)$ is therefore the groupoid $X \times X$
with the discrete topology.
I recall attending a lecture by Boris Schein in which he described the elements of $X \times X$ as being {\em infinitesimal elements}.
It was an attractive phrase but at the time I regarded it as metaphorical rather than mathematical.
I think the main theorem of this paper shows that it was in fact mathematical.}
\end{example}

\begin{example}{\em Our second example shows that our theory is analogous to the classical theory of sheaves of groups \cite{J}.
An inverse semigroup $S$ is said to be {\em Clifford} if its idempotents are central.
This is equivalent to the condition that $s^{-1}s = ss^{-1}$ for all $s \in S$.
Such inverse semigroups can be described as presheaves of groups over their semilattices of idempotents; see Section~5.2 of \cite{L} for details.
If $S$ is a boolean inverse monoid which is Clifford then the groupoid $\mathsf{G}(S)$
is a disjoint union of groups: to see why let $A$ be an ultrafilter in $S$.
We claim that $A^{-1} \cdot A = A \cdot A^{-1}$.
Let $x \in A^{-1} \cdot A$.
Then $a^{-1}a \leq x$ for some $a \in A$.
By assumption $a^{-1}a = aa^{-1}$.
Thus $aa^{-1} \leq x$ and so $x \in A \cdot A^{-1}$.
We have proved that $A^{-1} \cdot A \subseteq A \cdot A^{-1}$ and the reverse inclusion follows by symmetry.
It follows that in the groupoid $\mathsf{G}(S)$ the maps $\dom$ and $\ran$ coincide and we will write them both as $\mathbf{p}$.
We therefore have a local homeomorphism $\mathbf{p} \colon \mathsf{G}(S) \rightarrow \mathsf{G}(S)_{o}$ giving us a sheaf-space of groups.
We recover $S$ from this sheaf-space by using only those sections, which are in this context the same as bisections,
over the clopen subspaces of $\mathsf{G}(S)_{o}$.}
\end{example}

\begin{example} 
{\em The groupoids studied in \cite{H1,H2} are boolean when restricted to the case where the ultrametric spaces in question are compact.
From Lemma~3.5 \cite{H2}, the groupoids have a basis of open bisections each of which is homeomorphic to an open ball of the ultrametric space.
By Proposition~4.2(3) of \cite{H2}, every open ball is closed, and closed subsets of compact spaces are compact, thus each of these open bisections is also compact.
It follows that the groupoid has a basis of compact open bisections.}
\end{example}

\section{Cuntz groupoids}\setcounter{theorem}{0}
 
We begin by describing an example that connects our main theorem with work of ours on the 
way in which the Thompson groups $V_{n,1}$ may be constructed from the polycyclic inverse monoids \cite{L2,L3,L4}.
The polycyclic inverse monoids were first introduced and studied by Nivat and Perrot \cite{Nivat}.
They were rediscovered by Cuntz \cite{C} in the course of his work on what are now called Cuntz $C^{\ast}$-algebras;
for this reason, these inverse monoids are usually called {\em Cuntz inverse semigroups} in the $C^{\ast}$-algebra literature \cite{Paterson,Renault}.
Whatever one chooses to call them they are a fascinating class of inverse monoids arising both in formal language theory and the theory of wavelets \cite{BJ}.
We shall return to the polycyclic monoids later, but first we shall describe a class of boolean groupoids.

We shall use the following standard notation below.
If $A$ is a finite set then $A^{\ast}$ denotes the free monoid on $A$, which consists of all finite strings over $A$.
By $A^{\omega}$ we mean the set of all right-infinite strings over $A$.

We follow \cite{Paterson,Renault}.
For $n \geq 2$ and finite, put $A_{n} = \{a_{1}, \ldots, a_{n} \}$.
If $x$ is a finite string then $\left| x \right|$ is its length.
Then $G_{n}$ is the subset of $A_{n}^{\omega} \times \mathbb{Z} \times A_{n}^{\omega}$ consisting of triples
of the form $(xw, \left| x \right| - \left| y \right|, yw)$ where $x,y \in A_{n}^{\ast}$ and $w \in A_{n}^{\omega}$.
This set becomes a groupoid when the product is defined by $(z,k,z')(z',k',z'') = (z, k+k',z'')$ and inverses by $(z,k,z')^{-1} = (z',-k,z)$.
For $x,y \in A_{n}^{\ast}$ and $V \subseteq A_{n}^{\omega}$ open define
$$U_{x,y,V} = \{(xw, \left| x \right| - \left| y \right|, yw) \colon w \in V \},$$
clearly a bisection, and if $V = A_{n}^{\omega}$ denote $U_{x,y,A_{n}^{\omega}}$ by $U_{x,y}$.
Let $\Omega$ be the topology on $G_{n}$ with basis the sets $U_{x,y,V}$.
Then $G_{n}$ is a boolean groupoid called the {\em Cuntz groupoid}.
Our goal is to describe the boolean inverse monoid $\mathsf{A}(G_{n})$.
Let $P_{n}$ be the polycyclic monoid on $n$ generators; see \cite{L2} for a quick introduction.
Define $\psi \colon P_{n} \rightarrow \mathsf{A}(G_{n})$ by $\psi (xy^{-1}) = U_{x,y}$.
This is an injective homomorphism.
Each element of $\mathsf{A}(G_{n})$ is the finite orthogonal join of elements of $\psi (P_{n})$.
We now make one further observation.
The idempotents $a_{1}a_{1}^{-1}, \ldots, a_{n}a_{n}^{-1}$ form an orthogonal set in $P_{n}$.
Their images will have an orthogonal join in $\mathsf{A}(G)$ which we now calculate:
$$\bigvee_{i=1}^{n} \psi (a_{i}a_{i}^{-1}) 
= \bigvee_{i=1}^{n} U_{a_{i},a_{i}} 
= 1$$
the identity on the set of identities of $G_{n}$.
In the terminology of \cite{L3}, the inverse monoid $C_{n} = \mathsf{A}(G_{n})$ is the strong orthogonal
completion of the polycyclic monoid $P_{n}$.
We therefore have the following.

\begin{proposition} The boolean inverse monoid $C_{n}$ associated with the Cuntz groupoid $G_{n}$ is
the strong orthogonal completion of the polycyclic (or Cuntz) monoid $P_{n}$.
The group of units of $C_{n}$ is the Thompson group $V_{n,1}$.
\end{proposition}

To conclude this section, we carry out a calculation which suggests an avenue of further development.
All upwardly closecd cosets in $P_{n}$ were completely described in \cite{L4} motivated by calculations carried out by Kawamura \cite{Ka1,Ka2}.
We therefore have explicit descriptions of the ultrafilters in $P_{n}$.
In the result below, it is the groupoid, not the topological groupoid, which is constructed;
the topology is the one described in \cite{Lenz}. 

\begin{proposition} The elements of the groupoid $G_{n}$ can be identified with the ultrafilters in $P_{n}$.
\end{proposition}
\proof An idempotent ultrafilter $H$ in $P_{n}$ is determined by an element of $z \in A_{n}^{\omega}$
because
$$H = \{uu^{-1} \colon u \text{ is a finite prefix of } z \}^{\uparrow}.$$
Now let $A$ be an arbitrary ultrafilter such that $A^{-1} \cdot A = H$.
Then we may write $A = (a H)^{\uparrow}$ where $\dom (a) \in H$.
Let $a = xy^{-1}$ where $yy^{-1} \in H$ and so $y$ is a prefix of $z$.
Thus we may write $z = yw$.
We now calculate $K = A \cdot A^{-1}$.
This is just 
$$H = \{vv^{-1} \colon v \text{ is a finite prefix of } xw \}^{\uparrow}.$$
Thus the ultrafilter $A$ determines the ordered pair $(xw,yw)$.
However by choosing a different coset representative we obtain a different ordered pair.
Suppose that $x'y'^{-1}$ determines the same coset of $H$ as $xy^{-1}$.
Then $(x'y'^{-1})^{-1}xy^{-1} \in H$.
This product must be above an idempotent in $H$ and so, since $P_{n}^{\ast}$ is $E^{\ast}$-unitary, it must itself be a non-zero idempotent.
It follows that $x$ and $x'$ are prefix-comparable.
Assume that $x = x'p$ for some finite string $p$.
It follows that $(x'y'^{-1})^{-1}xy^{-1} = y'py^{-1}$.
Since this has to be an idempotent we have that $y'p = y$.
In addition, $y$ is a prefix of $z$.
Thus
$$(xw,\left| x \right| - \left| y \right|,yw) 
= 
(x'pw, \left| x' \right| - \left| y' \right|, y'pw).$$
In other words, the coset $A$ determines an element of the groupoid $G_{n}$.

Suppose now that
$$(\bar{z},k,z) = (xw,k,yw) = (x'w',k,y'w'),$$
an element of $G_{n}$.
Let $H$ be the idempotent ultrafilter of $P_{n}$ determined by the infinite string $z$.
We prove that $xy^{-1}$ and $x'y'^{-1}$ determine the same coset of $H$.
We shall suppose that $k \geq 0$.
Now $z = yw = y'w'$ and so $y$ and $y'$ are prefix-comparable.
Suppose that $y = y'p$ for some finite string $p$. 
Then $pw = w'$. 
Now $xw = x'w'$ and so $xw = x'pw$.
It is tempting to cancel the $w$.
But such a temptation must be resisted because $w$ is an infinite string.
It is here that we use the information provided by the number $k$.
We have that 
$$\left| x \right| 
= \left| y \right| + k 
= \left| y' \right| + \left| p \right| + k 
= \left| x' \right| + \left| p \right|
=  \left| x'p \right|.$$ 
But $x$ and $x'p$ are prefix-comparable and have the same length and so they must be equal.
We therefore have that $x = x'p$.
We can now calculate  
$(x'y'^{-1})^{-1}xy^{-1}$ which is a prefix of $z$ and so $xy^{-1}$ and $x'y'^{-1}$ determine the same coset of $H$.\qed \\

The significance of this example is that although $P_{n}$ is not a boolean inverse monoid it can still be used to construct a boolean groupoid.
From that groupoid one can construct a boolean inverse monoid, namely $C_{n}$, into which $P_{n}$ embeds.
The boolean inverse monoid $C_{n}$ is a completion of of $P_{n}$ as shown in \cite{L3}.
However, $C_{n}$ is somewhat complex.
To calculate the groupoid $G_{n}$ it is easier to start with $P_{n}$ and construct $G_{n}$ from the ultrafilters in $P_{n}$.
Such ultrafilters are intimately connected with what we call in \cite{LMS} {\em universal actions}.
We shall develop this idea further in a subsequent paper.


\end{document}